\documentclass[11 pt]{amsart}
\usepackage[all]{xy}
\usepackage{amssymb,amsmath,color}
\usepackage{amsfonts}
\usepackage[bottom]{footmisc}
\usepackage[normalem]{ulem}
\usepackage{amssymb}
\usepackage{amscd}
\usepackage{verbatim}
\usepackage{latexsym}
\usepackage{amsmath,color}
\usepackage{tikz-cd}
\usepackage{hyperref}

\newtheorem{theorem}{Theorem}

\newtheorem{lemma}[theorem]{Lemma}

\newtheorem{prop}[theorem]{Proposition}

\theoremstyle{remark}

\numberwithin{equation}{section}

\def\FF{{\mathbb F}}
\def\ff{{\mathbb F}}

\def\QQ{{\mathbb Q}}
\def\ZZ{{\mathbb Z}}

\def\ZZp{ { {\mathbb Z}/p{\mathbb Z}}}

\def\LP{ \left(\begin{array}}
\def\RP{ \end{array}\right)}
\def\LB{ \left[\begin{array}}
\def\RB{ \end{array}\right]}

\def\QQ{ {\mathbb Q}}

\def\Zp{ {\mathbb Z}_p }

\def\FF{ {\mathbb F} }

\def\Fp{ {\mathbb F}_p}

\def\Zp{ {\mathbb Z}/p{\mathbb Z} }

\DeclareFontFamily{U}{wncy}{}
\DeclareFontShape{U}{wncy}{m}{n}{<->wncyr10}{}
\DeclareSymbolFont{mcy}{U}{wncy}{m}{n}
\DeclareMathSymbol{\Sh}{\mathord}{mcy}{"58}

\catcode`,\active

\catcode`\,12

\title{On  tame $\ZZ/p\ZZ$-extensions with prescribed ramification}\thanks{
The second author  was partially supported by the ANR project FLAIR (ANR-17-CE40-0012) and  by the EIPHI Graduate School (ANR-17-EURE-0002). The third author was partially supported by Simons Collaboration grant \#524863. He also thanks FEMTO-ST for its hospitality and wonderful research environment during his visit there in the spring of 2022. All three authors were supported by ICERM for a Research in Pairs visit in January, 2022. }
\author{Farshid Hajir, Christian Maire, Ravi Ramakrishna}
\address{Department of Mathematics \& Statistics, University of Massachusetts, Amherst, MA 01003, USA}
 \address{FEMTO-ST Institute, Universit\'e Bourgogne Franche-Comt\'e, CNRS,  15B avenue des Montboucons, 25000 Besan\c con, FRANCE} 
\address{Department of Mathematics, Cornell University, Ithaca, NY 14853-4201USA}
\email{hajir@math.umass.edu, christian.maire@univ-fcomte.fr, ravi@math.cornell.edu}

\begin{document}
\maketitle\begin{abstract}
 %We give a new and \textcolor{black}{\sout{ elementary}} proof of 
The tame Gras-Munnier Theorem  gives a criterion for the existence of a $\Zp$-extension of a number field $K$ ramified at exactly a set $S$  of places of $K$ prime to $p$ (allowing real Archimedean places when $p=2$)  in terms of the existence of a dependence relation on the Frobenius elements of these places in a certain {\it governing extension}. 
%Using only one element of Class Field Theory and elementary linear algebra, 
We give a new and simpler proof of this theorem 
that also relates the set of such extensions of $K$ to the set of these dependence relations. After presenting this proof, we then reprove the key Proposition~\ref{prop:dims} using the  more sophisticated Wiles-Greenberg formula based on global duality. 
%as this argument may appear more natural.
\end{abstract}
\vspace{3mm}
%\noindent {\it Address:}

\noindent
\section{Introduction:} 
%\noindent Department of Mathematics, Cornell University, Malott Hall,
%Ithaca, NY 14853, USA. e-mail: {\tt ravi@math.cornell.edu}
Let $D\in \ZZ$ be  squarefree and odd. Our convention is that $\infty |D$ if $D<0$.
%\sout{It is an exercise in a beginning number theory class}
\textcolor{black}{It is a standard result} that there exists a quadratic extension
$K/\QQ$ ramified at exactly the set $\{v : v| D \}$ if and only if $D\equiv 1$ mod $4$. 
%If $D<0$ then we consider $\infty$ ramified in $K/\mathbb Q$.
%Of course if $p\equiv 3$ mod $4$ the quadratic extension $\QQ(\sqrt{-p})/\QQ$ has discriminant $-p$ and is ramified at both the one finite place $p$ {\it and} infinity. Let $\sigma_v$ denote Frobenius at $v$ in this extension with $\sigma_\infty$ being the nontrivial element of $Gal(\QQ(i)/\QQ)$.
The key is how the Frobenius elements of the $v|D$ lie in the Galois group of the  `governing extension' $\QQ(i)/\QQ$. 
Let $\sigma_v$ denote Frobenius at $v$ in this extension with $\sigma_\infty$ being the nontrivial element of $Gal(\QQ(i)/\QQ)$.
We frame this  as the following  theorem:  
%below is also an exercise.
%Proving the theorem below is an exercise:

\vskip1em\noindent
{\bf Theorem} {\it There exists a quadratic extension $K/\QQ$ ramified exactly at a tame (not containing $2$ 
\textcolor{black}{but allowing $\infty$}) set $S$ of places  if and only if $\displaystyle \sum_{v\in S} \sigma_v  $ is the trivial element in $Gal(\QQ(i)/\QQ)$.}

%\vskip1em
\noindent

In \cite{GM} this result was generalized to $\Zp$-extensions of a general number field $K$. 
For a fixed prime $p$ and set $S$ of tame places,  set
$$V_S = \{ x \in K^\times \mid (x)=J^p; \,\,x \in K^{\times p}_v \,\,\forall \,\,v \in S\}.$$
Note $K^{\times p} \subset V_S$ for all $S$ and $S\subset T \implies V_T \subseteq V_S$.
Let ${\mathcal O}^{\times}_K$ and $Cl_K[p]$ be, respectively, the units of $K$ and the $p$-torsion in the class group of $K$.
%Global class field theory implies  
\textcolor{black}{It is a standard result that } $V_\emptyset/K^{\times p}$ lies in the exact sequence (see Proposition 10.7.2 of \cite{NSW}):
$$0 \to  {\mathcal O}^{\times}_K \otimes \FF_p \to V_\emptyset / K^{\times p}  \to Cl_K[p] \to 0.$$
Set $K'=K(\mu_p)$,  $L=K'(\sqrt[p]{V_\emptyset})$ and let $r_1$ and $r_2$ be the number of real and pairs of complex embeddings of $K$.  We call $L/K'$ the {\it governing extension} for $K$.
When $K=\QQ$ and $p=2$ we see $L=\QQ(i)$ and have recovered the field of the theorem above. 
\begin{center}
\begin{tikzcd}
& L:=K'(\sqrt[p]{V_\emptyset}) \\
K':=K(\mu_p) \arrow[ur,dash]\\
K\arrow[u, dash]
\end{tikzcd}
\end{center}
%For any field $E$ set $\delta(E) = \left\{\begin{array}{cc} 1 & \mu_p \subset E\\0 &\mu_p \not \subset E \end{array}\right.$   and observe that by Dirichlet's unit theorem $Gal(L/K')$ is an $\FF_p$-vector space  of dimension $r_1+r_2-1+\delta(K)+ \dim (V_\emptyset/K^{\times p})$.  
%\vskip1em
%\hspace{2.3in}
%\begin{center}
%\begin{tikzcd}
%& L:=K'(\sqrt[p]{V_\emptyset}) \\
%K':=K(\mu_p) \arrow[ur,dash]\\
%K\arrow[u, dash]
%\end{tikzcd}
%\end{center}

As $L$ is obtained by adjoining to $K'$ the $p$th roots of elements of $K$ (not $K'$),   one easily shows that places  $v'_1,v'_2$ of $K'$ above a fixed place $v$ of $K$ have Frobenius elements in $ Gal(L/K')$ that differ by a nonzero scalar multiple.% \textcolor{black}{(necessarily $1$ when $p=2$)}.  %For the rest of this paper 
We abuse notation and for any $v'$ of $K'$ above $v$ in $K$ denote Frobenius at $v'$ by $\sigma_v$. The theorem of 
\cite{GM} \textcolor{black}{(also see Chapter V of \cite{gras})}
below and Theorem~\ref{theorem:main} implicitly use this abuse of notation.  %Our use of the term {\it place}  is only relevant when $p=2$ and $K_v=\mathbb R$  and by Frobenius of such a $v$ we mean complex conjugation in $Gal(L_w/K_v)=Gal({\mathbb C}/{\mathbb R})$. Frobenius at $v$ is denoted by $\sigma_v$ in the Archimedean case as well.

\vskip1em\noindent
{\bf Theorem} (Gras-Munnier) {\it Let $p$ be a prime and $S$ a finite  set of {\it tame} places \textcolor{black}{ (allowing Archimedean places if $p=2$)} of $K$. There exists a $\ZZ/p\ZZ$-extension of $K$ ramified at {\it exactly} the places of $S$ if and only if  there exists a dependence relation 
$\displaystyle \sum_{v\in S} a_v \sigma_v =0$ in the $\FF_p$-vector space $Gal(L/K')$ with all $a_v \neq 0$.}

%\vskip1em
\noindent

The original proof  uses class field theory in a fairly complicated way. Theorem~\ref{theorem:main} 
is a  generalization of the Gras-Munnier Theorem. We first give a short proof 
that uses only  one element of class field theory, \eqref{eq:cft} below,
and elementary linear algebra.
% on the space of $\ZZ/p\ZZ$ extensions of a number field ramified unramified outside a tame set $Z$.
%The trade-off is that the rest of our proof involves elementary, but moderately complicated linear algebra over $\FF_p$. 
We easily prove Proposition~\ref{prop:dims} from \eqref{eq:cft}, after which one only needs a standard inclusion-exclusion argument. 
\textcolor{black}{The cardinalities of the two sets of Theorem~\ref{theorem:main} being equal suggests a duality. 
In the final section of this note we give an alternative proof of Proposition~\ref{prop:dims} using the  Wiles-Greenberg formula whose proof requires the full strength of global duality.}
\textcolor{black}{Denote by $G_S$ the Galois group over $K$ of its maximal extension pro-$p$ unramified outside $S$.}
%As we will see, there is one striking difference between the $p=2$ and $p>2$ cases.
Our main result is:

\begin{theorem}\label{theorem:main}
Let $p$ be a prime and $S$  a finite set of {\it tame} places of a number field $K$ \textcolor{black}{ (allowing Archimedean places if $p=2$)}. The sets
$$\left\{ f\in \frac{H^1(G_S,\Zp)}{H^1(G_\emptyset,\Zp)} \,\mid  \, \mbox{ the extension } K_f/K \mbox{ fixed by } Ker(f) \mbox{ is ramified exactly at the places of } S\right\}$$
and 
$$\{\mbox{The dependence relations }\sum_{v\in S} a_v\sigma_v=0 \mbox{ in } Gal(L/K')\mbox{ with all } a_v\neq 0 \}$$
have the same cardinality.
%For $p=2$ there is an obvious natural bijection between these sets.
\end{theorem}
It is an easy exercise to see both sets have cardinality at most one when $p=2$, so the bijection is natural in this case.
%\sout{Of course $\ff^{\times}_p$ acts on both  sets above  so we can quotient out by that action. When $p=2$ there is no scaling and the bijection is clear.}

\textcolor{black}{We thank Brian Conrad for pointing out to us a proof of Lemma~\ref{lemma:XXX}(ii) and Peter Uttenthal for helpful suggestions.}

\section{Proof of the Gras-Munnier Theorem}
%\vskip1em
\textcolor{black}{For any field $E$ set $\delta(E) = \left\{\begin{array}{cc} 1 & \mu_p \subset E\\0 &\mu_p \not \subset E \end{array}\right. .$   Observe that by Dirichlet's unit theorem $Gal(L/K')$ is an $\FF_p$-vector space  of dimension $r_1+r_2-1+\delta(K)+ \dim Cl_K[p]$. }
The standard fact from class field theory that we need  (see \S 11.3 of \cite{K} or \S 10.7 of \cite{NSW}) 
\textcolor{black}{is a formula of Shafarevich and Koch} for 
the dimension of the space of  $\Zp$-extensions of $K$  unramified outside a tame set $Z$:
\begin{equation}\label{eq:cft}
 \dim H^1(G_Z,\Zp) =-r_1-r_2+1-\delta(K) + \dim (V_Z/K^{\times p}) + \left(\sum_{v \in Z} \delta(K_v) \right) .
\end{equation}

Fix a tame set $S$ noting that $H^1(G_S,\Zp)$ may include cohomology classes that cut out 
$\Zp$-extensions of $K$ that could be ramified at proper subsets of $S$. 
If $\delta(K_v) =0$ for $v \in S$, there are no ramified $\Zp$-extensions of $K_v$ and thus no $\Zp$-extensions of $K$ ramified at $v$, so we always assume $\delta(K_v)=1$. Then, as we vary  $Z$  from $\emptyset$ to $S$ one place at a time, $\dim (V_Z/K^{\times p})$ may remain the same or decrease by $1$. In these cases $\dim H^1(G_Z,\Zp)$ increases by $1$ or remains the same respectively.

 Let $W \subset Gal(L/K')$ be the $\FF_p$-subspace spanned by $\langle \sigma_v\rangle_{v\in S}$,
 \textcolor{black}{the Frobenius elements of the places in $S$. We will show the set of dependence relations on these Frobenius elements all of whose coefficients are nonzero has the same cardinality as the set of $\ZZ/p\ZZ$-extensions of $K$ ramified exactly at the places of $S$.}
 Let 
% $I:=\{ \p_1, \p_2,\cdots,\p_r\} \subset S $ 
 $I:= \{u_1,u_2,\cdots,u_r\} \subset S$
 be such that  
 $\{\sigma_{u_1},\sigma_{u_2},\cdots,\sigma_{u_r} \}$ 
%$\{u_1,u_2,\cdots,u_r\}$
 form a  basis of $W$ and let 
 %$D:=\{ \q_1,\q_2,\cdots ,\q_s\} \subset S$ 
 $D:=\{ w_1,w_2,\cdots ,w_s\} \subset S$
 be the  remaining elements of $S$. We think of the $\sigma_{u_i}$ as independent elements
and the $\sigma_{w_j}$ as the dependent elements.
As we vary $X$ in \eqref{eq:cft} from $\emptyset$ to $I$ by adding in one $u_i$ at a time, we are adding $1$ through the $\delta(K_{u_i})$ term
to the right side of \eqref{eq:cft},
but  $\dim V_X/K^{\times p}$ becomes
one dimension smaller, so both sides remain unchanged. 
Then, as we add in the places $w_j$ of $D$ to get to $S=I\cup D$ we have $V_I/K^{\times p} = V_S/K^{\times p}$.
Thus 
\begin{equation}\label{eq:h1eq}
H^1(G_\emptyset,\Zp) = H^1(G_I,\Zp)=\dim H^1(G_S,\Zp) - s \implies \dim \frac{H^1(G_S,\Zp)}{H^1(G_\emptyset,\Zp)}=s.
\end{equation}

We can write each $\sigma_{w_j}$ as a linear combination of the 
$\sigma_{u_i}$
in a unique way: 
%\begin{equation}\label{depend}
$$R_j:\,\, \sigma_{w_j} - \sum^r_{i=1} F_{ji} \sigma_{u_i}=0. $$
%\mbox{ or equivalently }
%\left[\begin{array}{c} \sigma_{\q_1}\\ \sigma_{\q_2} \\ \vdots \\ \sigma_{\q_s} \end{array}\right]= [F_{jk}]
%\left[\begin{array}{c} \sigma_{\p_1}\\ \sigma_{\p_2} \\ \vdots \\ \sigma_{\p_r} \end{array}\right]
%\end{equation}
%Let $S_j=\{\q_j\} \cup \{ \p_k | F_{jk}\neq 0\}$, the support of the dependence relation $R_j$. 
For $X\subset S$ let $R_X$ be the $\FF_p$-vector space of all dependence relations on the elements 
$\{\sigma_v\}_{v\in X} \subset Gal(L/K')$.
We prove a preliminary result  in the spirit of Theorem~\ref{theorem:main}.
\begin{lemma}\label{lemma:dep}
The set $\{R_1,R_2,\cdots,R_s\}$ forms a basis of the \textcolor{black}{$\FF_p$-vector space of} dependence relations on the $\sigma_{u_i}$ and $\sigma_{w_j}$.\\
\end{lemma}
\begin{proof}
Consider {\it any} dependence relation $R$ among the $\sigma_{u_i}$ and $\sigma_{w_j}$. We can eliminate each $\sigma_{w_j}$ by adding to $R$ a suitable multiple of $R_j$.  We are then left with a dependence relation on the $\sigma_{u_i}$, which are independent, so it is trivial, proving the lemma.
%Thus $R$ can be uniquely written as a linear combination of the $R_j$.
\end{proof}

\begin{prop} \label{prop:dims} For any $X \subseteq S$,
%Let $R_X$ be the space of dependence relations with support contained in $X$. Then 
$\displaystyle\dim R_X = \dim \left( \frac{H^1(G_X,\Zp)}{H^1(G_\emptyset,\Zp)}\right).$
\end{prop}
\begin{proof} Lemma~\ref{lemma:dep} and ~\eqref{eq:h1eq} prove this for $X=S$. Apply the same proof to $X\subset S$.
\end{proof}

\vskip1em\noindent
 \textcolor{black}{\it Proof of Theorem~\ref{theorem:main}}:
%\begin{remark}
Proposition~\ref{prop:dims} does {\it not} complete the proof  as $R_S$ may contain dependence relations with support properly contained in $S$ and
$\displaystyle \frac{H^1(G_S,\Zp)}{H^1(G_\emptyset,\Zp)}$ may contain elements giving rise to extensions of $K$ ramified at proper subsets of $S$. 
%\end{remark}
\vskip1em\noindent
{\it Proof Theorem~\ref{theorem:main}:}
The set of dependence relations with support {\it exactly} in $S$ is 
\begin{equation}\label{eq:rs}
R_S \backslash \bigcup_{v \in S} R_{S \backslash \{v\}},
\end{equation}
 those with support contained in $S$ less \textcolor{black}{the union of} those with proper maximal support  in $S$. 
%These last relations are those whose support is $S$ less one prime. 
%The cardinality of the union of these latter subspaces is determined 
For any sets $A_i \subset S$ \textcolor{black}{it is clear that} 
%\begin{equation}\label{eq:rai}
$$\bigcap^k_{i=1} R_{A_i} = R_{\bigcap^k_{i=1} A_i},$$
%\end{equation}
so
by inclusion-exclusion 
\begin{equation}\label{eq:rai2}  \# \bigcup_{v \in S} R_{S \backslash \{v\}} =\sum_{v\in S} \# R_{S \backslash \{v\}} - 
\sum_{v\neq w \in S} \# R_{S \backslash \{v,w\}} +\cdots
\end{equation}
%and is an alternating sum of the cardinalities of their intersections.

Similarly the set of cohomology classes giving rise to $\ZZ/p\ZZ$-extensions ramified exactly at the places of $S$ (up to unramified extensions) is 
\begin{equation}\label{eq:aai2}
\frac{H^1(G_{S},\Zp)}{H^1(G_\emptyset,\Zp)}  \backslash \bigcup_{v\in S} \frac{ H^1(G_{ S\backslash \{v\} },\Zp)}{H^1(G_\emptyset,\Zp)}.
\end{equation}
%in bijection with the linear combinations of the $f_i$ that are unramified outside $S$ less those that are unramified outside proper maximal subsets of $S$.  Again, the cardinality of the union of these latter subspaces is determined by inclusion-exclusion and 
%depends only on the cardinalities of the these subspaces and their various intersections. 
Since for any sets $A_i \subset S$ we have
%\begin{equation}\label{eq:aai}
$$\bigcap^k_{i=1} \frac{H^1(G_{A_i},\Zp)}{H^1(G_\emptyset,\Zp)} =\frac{ H^1(G_{ \cap^k_{i=1} A_i},\Zp)}{H^1(G_\emptyset,\Zp)},$$
%\end{equation}
we see 
\begin{equation}\label{eq:aai3}
\# \bigcup_{v\in S} \frac{ H^1(G_{ S\backslash \{v\} },\Zp)}{H^1(G_\emptyset,\Zp)} = \sum_{v\in S} 
\# \frac{ H^1(G_{ S\backslash \{v\} },\Zp)}{H^1(G_\emptyset,\Zp)} 
- \sum_{v\neq w \in S} \# \frac{ H^1(G_{ S\backslash \{v,w\} },\Zp)}{H^1(G_\emptyset,\Zp)}  + \cdots
\end{equation}
Proposition~\ref{prop:dims}  implies the terms on the right sides of  \eqref{eq:rai2} and \eqref{eq:aai3} are equal so the left sides are equal as well. The result follows from~\eqref{eq:rs}, ~\eqref{eq:aai2} and applying Proposition~\ref{prop:dims} with $X=S$.
%imply each $k$-fold intersection of spaces of dependence relations with maximal proper support has the same dimension as the corresponding intersection of cohomology classes with maximal proper ramification  (up to unramified classed) and hence these intersections have the same cardinalities.  The result follows from equations ~\eqref{eq:rs} and ~\eqref{eq:aai2}.
%By inclusion-exclusion  the set consisting of the union of dependence  relations with  proper maximal  support has the same cardinality as the union of cohomology classes with proper maximal ramification. These unions of equal cardinality are subtracted from the equal cardinality sets $R_S$ and $H^1(G_S,\Zp)$ (Proposition~\ref{prop:dims}).  Thus the set of dependence relations with support exactly $S$  has the same cardinality as those $f\in H^1(G_S,\Zp)$ giving rise to extensions ramified exactly at $S$ (up to unramified classes).  Theorem~\ref{theorem:main} is proved.
%This proves tame Gras-Munnier.
 \hfill $\square$
 
\section{A proof via the Wiles-Greenberg formula}
  
%  \vskip1em\noindent\underline{A proof via the Wiles-Greenberg formula}: 
As the association of dependence relations and cohomology classes in Theorem~\ref{theorem:main} resembles a duality result, we now prove Proposition~\ref{prop:dims} using the Wiles-Greenberg formula, which follows from global duality. 
%The rest of the proof of Theorem~\ref{theorem:main} uses the inclusion-exclusion argument above. 
We assume familiarity with local and global Galois cohomology and their duality theories.

\textcolor{black}{As we will need to apply the Wiles-Greenberg formula, we henceforth assume its hypothesis that
$Z$} is a set of places of $K$ containing all those above infinity and $p$. 
%\sout{and $G_Z$ be the Galois group of the maximal extension of $K$ unramified outside $Z$}. 
For each $v\in Z$, let $G_v:=Gal(\bar{K}_v/K_v)$ and consider a subspace $L_v \subset H^1(G_v,\ZZ/p\ZZ)$.  Under the local duality pairing (see Chapter 7, \S 2 of \cite{NSW})
 $$H^1(G_v,\ZZ/p\ZZ) \times H^1(G_v,\mu_p) \to H^2(G_v,\mu_p) \simeq \frac1p\ZZ/\ZZ$$
$L_v$ has an annihilator $L^\perp_v \subset H^1(G_v,\mu_p)$.
Set $$H^1_{\mathcal L}(G_Z,\ZZ / p\ZZ) := \mbox{Ker} \left( H^1(G_Z,\ZZ / p\ZZ) \to \oplus_{v\in Z} 
\frac{H^1(G_v,\ZZp)}{L_v}\right)$$
and 
$$H^1_{{\mathcal L}^\perp}(G_Z,\mu_p) := \mbox{Ker} \left( H^1(G_Z,\mu_p) \to \oplus_{v\in Z} 
\frac{H^1(G_v,\mu_p)}{L^\perp_v}\right).$$
We call $\{L_v\}_{v \in Z}$ and $\{ L^\perp_v \}_{v\in Z}$  the Selmer and   dual Selmer conditions and 
$H^1_{\mathcal L}(G_Z,\ZZ / p\ZZ)$ and $H^1_{{\mathcal L}^\perp}(G_Z,\mu_p)$ the Selmer and dual Selmer groups.

We state two results that we need for our second proof of Proposition~\ref{prop:dims}. As  Lemma~\ref{lemma:XXX} (ii) is perhaps not so well-known, we include \textcolor{black}{ a sketch of its proof.}
%For $v$ tame or infinite, we will set $L_v=H^1(G_v,\ZZ/p\ZZ)$ (so $L^\perp_v = 0$) as we ramification restrictions at such $v$ for  those elements of $H^1(G_S,\ZZ/p\ZZ)$  we wish to capture in our Selmer group .  For $v|p$ we set $L_v=H^1_{nr}(G_v,\ZZ/p\ZZ)$ (so $L^\perp_v = H^1_f(G_v,\ZZ/p\ZZ)$, see. Lemma~\ref{lemma:XXX} below) as we only wish to capture classes that are unramified at these places.

%It is well-know that for $v\nmid p$ that the annihilator of $H_{nr}^1(G_v,\ZZ/p\ZZ)$ is $H_{nr}^1(G_v,\mu_p)$.
%The situation for $v|p$ is somewhat different.

%Part (i) of Lemma~\ref{lemma:XXX} below is a standard result and we do not give the proof. Part (ii)  seems to be known to the experts but we have found no proof in the literature.
%We thank Brian Conrad for pointing out the proof below to us.

\begin{lemma} \label{lemma:XXX} 
(i) Suppose $v\nmid p$. Then $H_{nr}^1(G_v,\ZZ/p\ZZ)$ and $H_{nr}^1(G_v,\mu_p)$, the unramified cohomology classes, are exact annihilators of one another under the local duality pairing.\\
(ii) Suppose $v|p$. \textcolor{black}{Recall $K'_v=K_v(\mu_p)$.} The annihilator in 
$H^1(G_v,\mu_p)$ of $H_{nr}^1(G_v,\ZZ/p\ZZ) \subset H^1(G_v,\ZZ/p\ZZ)$ is $H^1_f(G_v,\mu_p)$, the peu ramifi\'{e}e classes, namely those $f\in H^1_f(G_v,\mu_p)$ whose fixed field $L_{v,f}$ of Kernel$(f |_{G_{K'_v}})$ arises from adjoining the $p$th root of a unit $u_f \in K_v$.
%
%\vskip1em
%\hspace{2.3in}
%\begin{tikzcd}
%& L_{v,f}:=K_v'(\sqrt[p]{u_f}) \\
%K_v':=K_v(\mu_p) \arrow[ur,dash]\\
%K_v\arrow[u, dash]
%\end{tikzcd}
%\vskip1em
\end{lemma} 

\begin{proof}
(i) This is standard - see 7.2.15 of \cite{NSW}. \newline\noindent
(ii) %We thank Brian Conrad for pointing out this proof to us.
 \textcolor{black}{Cohomology taken over $Spec({\mathcal O}_{K_v})${ in what follows is flat.}
Here $$H^1_{f}(G_v,\mu_p) =
H^1(Spec({\mathcal O}_{K_v}),\mu_p) = { {\mathcal O}^\times_{K_v}}/  { {\mathcal O}^{\times p}_{K_v}}
 \subset K^{\times}_v/K^{\times p}_v$$
where the containment is codimension one as $\Fp$-vector spaces. Recall $$\ZZp \simeq H^1_{nr}(G_v,\ZZp) 
= H^1(Spec({\mathcal O}_{K_v}),\ZZp)$$  and by Lemma 1.1 of Chapter III of \cite{Milne} we have the injections
$$H^1(Spec({\mathcal O}_{K_v},\ZZp)) \hookrightarrow H^1(G_v,\ZZp) \mbox{ and }   H^1(Spec({\mathcal O}_{K_v},\mu_p)) \hookrightarrow H^1(G_v,\mu_p)$$
%By Lemma 1.1 of Chapter III of \cite{Milne} we have 
and the pairing 
$$H^1(Spec({\mathcal O}_{K_v}),\ZZp) \times H^1(Spec({\mathcal O}_{K_v}),\mu_p) \to H^2(Spec({\mathcal O}_{K_v}),\mu_p)=0$$ 
which is consistent with the local duality pairing
$$H^1(G_v,\ZZp)\times H^1(G_v,\mu_p) \to H^2(G_v,\mu_p)=\frac1{p}\ZZ /\ZZ.$$
As $H^1(Spec({\mathcal O}_{K_v}),\ZZp)= H^1_{nr}(G_v,\ZZ/p\ZZ)$ and $H^1(Spec({\mathcal O}_{K_v}),\mu_p) =H^1_f(G_v,\mu_p)$ are, respectively,  dimension $1$ and codimension $1$ in
$H^1(G_v,\ZZp)$ and $H^1(G_v,\mu_p)$, they are exact annihilators of one another,
%each group is the exact annihilator of the other  under the local duality pairing $$H^1(G_v,\ZZp)\times H^1(G_v,\mu_p) \to H^2(G_v,\mu_p)=\ZZp,$$ 
proving (ii).}
\end{proof}

\noindent
{\bf Theorem} {\it (Wiles-Greenberg) \textcolor{black}{Assume $Z$ contains all places above $\{p,\infty\}$. Then} 
%\sout{With $Z$ as above,}
%Assume $Z$ contains all primes of $K$ above $p$ and all Archimedean places if $p=2$. Then
\begin{eqnarray*} 
\begin{tabular}{ll}
$\dim H^1_{\mathcal L}(G_Z,\ZZ / p\ZZ) - \dim H^1_{{\mathcal L}^\perp}(G_Z,\mu_p)$ & \\
$=\dim H^0(G_Z,\ZZ/p\ZZ) - \dim H^0(G_Z,\mu_p) + \sum_{v\in Z} \left(\dim L_v - \dim H^0(G_v,\ZZ/p\ZZ) \right).$
&
\end{tabular}
\end{eqnarray*}}
See $8.7.9$ of \cite{NSW} for details of this result.
%We will use the Wiles-Greenberg formula to give an alternative proof to Proposition~\ref{prop:dims}.
\vskip1em
\noindent
{\it Second proof of Proposition~\ref{prop:dims}.}
%We know $\dim L^\perp_v +1=\dim H^1(G_v,\mu_p)$. 
%Recall that for $v\nmid p$, we always assume $\mu_p \subset K_v$.
Recall $X$ is tame and write $X:=X_{<\infty} \cup X_\infty$. %The Archimedean places only play a role when $p=2$. 
Set $Z:= Z_p \cup X_{<\infty} \cup Z_\infty$ where
 $Z_p:= \{v : v|p\}$ and $Z_\infty$ is the set of all real Archimedean places of $K$ (so $X_\infty \subseteq Z_\infty$).
 % and $X_{<\infty}:=X_{<\infty}$. 
We assume  for all $v\in X_{<\infty}$ that $N(v) \equiv 1\mod p$.
%Put $\tilde{X}=X\cup X_p \cup X_\infty$.
%We define two sets of Selmer conditions, $\{M_v\}_{v\in Z}$ and $\{N_v\}_{v\in Z}$ so that $H^1_{\mathcal M}(G_{Z},\ZZp) =H^1(G_\emptyset,\ZZp)$ and  $H^1_{\mathcal N}(G_{Z},\ZZp) =H^1(G_{X},\ZZp)$. 

Recall that for a complex Archimedean place $v$ of $K$ we have $G_v=\{e\}$ so the Selmer and dual Selmer conditions are trivial. For a real Archimedean place $v$, $\dim H^1(G_v,\ZZ/2\ZZ) = \dim H^1(G_v,\mu_2)=1$ and the pairing between them is perfect - see Chapter I, Theorem $2.13$ of \cite{Milne}. It is easy to see in this case that the unramified cohomology groups are trivial. 

In the table below we choose $\{M_v\}_{v\in Z}$ and $\{N_v\}_{v\in Z}$ so that
$H^1_{\mathcal M}(G_{Z},\ZZp) =H^1(G_\emptyset,\ZZp)$ and 
$H^1_{\mathcal N}(G_{Z},\ZZp) =H^1(G_{X},\ZZp)$. \textcolor{black}{The previous paragraph and  Lemma~\ref{lemma:XXX}
justify
the stated dual Selmer conditions of the table.}
\begin{eqnarray*} 
\begin{tabular}{l|cccc}
& $M_v$ & $M^\perp_v$ & $N_v$ & $N^\perp_v$ \\
\hline
$v\in Z_p$ & $H^1_{nr}(G_v,\ZZ/p\ZZ)$ & $H^1_{f}(G_v,\mu_p)$  & $H^1_{nr}(G_v,\ZZ/p\ZZ)$ & $H^1_{f}(G_v,\mu_p)$\\
$v\in X_\infty$ & $H^1_{nr}(G_v,\ZZ/2\ZZ)=0$ & $H^1(G_v,\mu_2)$ & $H^1(G_v,\ZZ/2\ZZ)$ &  $0$ \\
$v\in Z_\infty \backslash X_\infty$ & $H^1_{nr}(G_v,\ZZ/2\ZZ) =0$ & $H^1(G_v,\mu_2)$ & $0$ & $H^1(G_v,\mu_2)$ \\
$v\in X_{<\infty}$ &  $H^1_{nr}(G_v,\ZZ/p\ZZ)$ & $ H^1_{nr}(G_v,\mu_p)$ & $H^1(G_v,\ZZ/p\ZZ)$ & $0$\\
\end{tabular}
\end{eqnarray*}
%\vskip1em
%\begin{itemize}
%\item $L_v$: \begin{itemize} 
%                         \item For $v \in Z_p$ set $L_v=H^1_{nr}(G_v,\ZZ/p\ZZ)$ so $L^\perp_v=H^1_{f}(G_v,\mu_p)$.
%                         \item For $v \in Z_\infty$ set $L_v=0$ so $L^\perp_v=H^1(G_v,\ZZ/p\ZZ)$.
%                         \item For $v \in X_{<\infty}$ set $L_v=H^1_{nr}(G_v,\ZZ/p\ZZ)$ so $L^\perp_v=H^1_{nr}(G_v,\mu_p)$.
%                     \end{itemize}
%\item $N_v$: \begin{itemize} 
%                        
%                         
%                         \item For $v \in Z_p$ set $N_v=H^1_{nr}(G_v,\ZZp)$ so $N^\perp_v=H^1_{f}(G_v,\mu_p)$.
%                          \item For $v \in Z_\infty$ set $L_v=H^1(G_v,\ZZ/p\ZZ)$ so $L^\perp_v=H^1_{nr}(G_v,\mu_p)$.
%                          \item For $v \in X_{<\infty}$ set $N_v=H^1(G_v,\ZZ/p\ZZ)$ so $N^\perp_v=0$.
%                     \end{itemize}
%
%
%\end{itemize}
%
%The point of going through these machinations is that we can now use use the Wiles-Greenberg formula, a result that comes from global duality and requires wild ramification.
%For all $v \in Z_p \cup X_{<\infty}$ we easily see $\dim H_{nr}^1(G_v,\ZZp) = 1$ so its annihilator in $H^1(G_v,\mu_p)$ has codimension  $1$ (except when $v|\infty$ and $p>2$ in which case the groups are trivial).     
%Thus for $v|p$ we have $\dim L^\perp_v = [K_v:\QQ_p] +\delta(K_v)$ and for $v\nmid p$ we have $\dim L^\perp_v = 1$.

Applying the Wiles-Greenberg formula for $\{M_v\}_{v\in Z}$ and $\{N_v\}_{v\in Z}$  and subtracting the first equation from the second:
\begin{eqnarray*}
\begin{tabular}{l}
$\dim H^1(G_X,\ZZp) -\dim H^1(G_\emptyset,\ZZp)=  $\\
$\dim H^1_{\mathcal N}(G_{Z},\ZZp) - \dim H^1_{\mathcal M}(G_{Z},\ZZp) =$\\
$\dim H^1_{{\mathcal N}^\perp}(G_{Z},\mu_p) - \dim H^1_{{\mathcal M}^\perp}(G_{Z},\mu_p) 
+\sum_{v\in Z} (\dim N_v - \dim M_v)$
%$\dim H^1_{{\mathcal N}^\perp}(G_{Z},\mu_p) - \dim H^1_{{\mathcal L}^\perp}(G_{Z},\mu_p) +\#(Z_\infty \cup X_{<\infty}).$
\end{tabular}
\end{eqnarray*}
%By local class field theory, 
For $v \in X_{<\infty}$ local class field theory implies  $\dim H^1_{nr}(G_v,\ZZp) =1$ and $\dim H^1(G_v,\ZZp) =2$ so
$$\dim N_v - \dim M_v =\left\{ 
\begin{array}{ll} 
0 & v\in Z_p \\
%0 & v\in Z_\infty \backslash X_\infty \\
%0 & v \in X_\infty, p\neq 2\\
1 & v \in X_\infty, \mbox{ }p=2\\
0 & v\in Z_\infty \backslash X_\infty \\
1 & v\in  X_{<\infty}
\end{array} \right. ,$$ 
and we have
\begin{equation}\label{eq:prop32}
\dim \left(\frac{H^1(G_X,\ZZp)}{H^1(G_\emptyset,\ZZp)} \right)= \dim H^1_{{\mathcal N}^\perp}(G_{Z},\mu_p) - \dim H^1_{{\mathcal M}^\perp}(G_{Z},\mu_p) +\#X .
\end{equation}
%For $v\in Z_p$ the $\dim N_v - \dim L_v$ terms are $0$. For $v\in Z_\infty$ they are $1$ and using the local Euler-Poincar\'{e} characteristic, we see they are $1$ for $v\in X_{<\infty}$ as well. Thus:

%Set $K':=K(\mu_p)$. 
To prove Proposition~\ref{prop:dims} we need to show this last quantity is $\dim R_X=s$, the dimension of the space of dependence relations on the set $\{\sigma_v\}_{v\in X} \subset W = 
Gal(K'(\sqrt[p]{V_\emptyset})/K')$.

An element $f \in H^1_{{\mathcal M}^\perp}(G_{Z},\mu_p)$ gives rise to the field diagram below
where $L_f/K'$ is a $\ZZp$ extension  peu ramifi\'{e}e at $v\in Z_p$, with no condition on $v\in Z_\infty$ and unramified at $v\in X_{<\infty}$.
We show the composite of all such $L_f$ is $K'\left(\sqrt[p]{V_\emptyset}\right)$.

\vskip1em
\hspace{2.3in}
\begin{tikzcd}
& L_{f}:=K'\left(\sqrt[p]{\alpha_f}\right) \\
K':=K(\mu_p) \arrow[ur,dash]\\
K\arrow[u, dash]
\end{tikzcd}

\vskip1em\noindent
\textcolor{black}{Kummer Theory  implies $\alpha_f \in K'/{K'}^{\times p}$, which decomposes into eigenspaces under the action of
$Gal(K'/K)$. If it is not in the trivial eigenspace, 
%We show $\alpha_f \in K$. $\alpha_f \in K' \backslash K$  
then $Gal(L_f/K')$ is not acted  on by $Gal(K'/K)$ via the cyclotomic character, a contradiction, so 
we may assume (up to $p$th powers) $\alpha_f \in K$.}
 Since $L_f/K'$ is unramified at $v\in X_{<\infty}$, we see that at all such $v$ that $\alpha_f = u \pi_v^{pr}$ where $u\in K_v$ is a unit. At $v\in Z_p$ being peu ramifi\'{e}e implies  that locally at 
$v\in X_p$ we have $\alpha_f = u \pi_v^{pr}$ where $u\in K_v$ is a unit. Together, these mean that the fractional ideal $(\alpha_f)$ of $K$ is a $p$th power, which implies that $\alpha_f \in V_\emptyset$. Conversely, if $\alpha \in V_\emptyset$, then, recalling that
$(\alpha)=J^p$ for some ideal of $K$, we have that
$K'\left(\sqrt[p]{\alpha}\right)/K'$ is a $\ZZp$-extension peu ramifi\'{e}e at $v\in Z_p$, with no condition at $v\in Z_\infty$. 
Thus $\alpha $ gives rise to an element $f_\alpha \in H^1_{{\mathcal M}^\perp}(G_Z,\mu_p)$
%and unramified at $v\in X_{<\infty}$, 
so $L:=K'\left( \sqrt[p]{V_\emptyset}\right)$ is the composite of all $L_f$ for 
$f \in H^1_{{\mathcal M}^\perp}(G_Z,\mu_p)$.

An element $f \in H^1_{{\mathcal N}^\perp}(G_{Z},\mu_p)$ gives rise to  a $\ZZp$-extension of $K'$
peu ramifi\'{e}e at $v\in Z_p$ and split completely at $v\in X$. We denote the composite of all these fields by
$D \subset K'\left( \sqrt[p]{V_\emptyset}\right)$.

\vskip1em
\hspace{2.3in}
\begin{tikzcd}
& & L:=K'\left( \sqrt[p]{V_\emptyset}\right) \\
& D\arrow[ur,dash] \\
K':=K(\mu_p) \arrow[ur,dash]\\
K\arrow[u, dash]
\end{tikzcd}
\vskip1em
\noindent
Recall that $r$ is the dimension of the space $\langle \sigma_v \rangle_{v\in X} \subset Gal(L/K')$.
Clearly $D$ is the field fixed of $\langle \sigma_v\rangle_{v\in X}$ so
$\dim_{\ff_p} Gal\left( K'\left( \sqrt[p]{V_\emptyset}\right)/ D \right) =r =\# I$ from the first section of this note.
Thus $$\dim H^1_{{\mathcal N}^\perp}(G_{Z},\mu_p)=\dim (V_\emptyset/K^{\times p}) -r $$
so 
\begin{eqnarray*}
\dim H^1_{{\mathcal N}^\perp}(G_{Z},\mu_p) - \dim H^1_{{\mathcal M}^\perp}(G_{Z},\mu_p) 
+\#X=\\
\left(\dim (V_\emptyset/K^{\times p}) -r \right) - \dim (V_\emptyset/K^{\times p}) +(r+s) =s = \dim R_X
\end{eqnarray*}
and we have shown the the left hand side of~\eqref{eq:prop32} is $\dim R_X$ proving Proposition~\ref{prop:dims}.
\hfill$\square$
\vskip1em
We have now proven Proposition~\ref{prop:dims} using the Wiles-Greenberg formula. The rest of the proof of Theorem~\ref{theorem:main} follows as in the previous section.

\end{document}